\font\medbf=cmbx10 scaled \magstep1
\font\bigbf=cmbx10 scaled \magstep2
\magnification=1000

\input psfig.sty
\hfuzz=3 pt
\hoffset=0.1in
\tolerance 9000
\def\i{\item}

\def\v{\vskip 1em}
\def\vs{\vskip 2em}
\def\vsk{\vskip 4em}
\def\c{\centerline}
\def\a{\alpha}
\def\i{\item}

\def\ov{\overline}

\def\bfw{{\bf w}}

\def\R{\Bbb R}

\def\A{{\cal A}}

\def\Om{\Omega}
\def\C{{\cal C}}
\def\D{{\cal D}}

\def\C{{\cal C}}
\def\H{{\cal H}}

\def\bfU{{\bf U}}

\def\R{I\!\!R}
\def\forall{\hbox{for all~~}}

\def\bfn{\hbox{\bf n}}
\def\ve{\varepsilon}
\def\n{\noindent}

\def\rhead=\vbox to 2truecm{\hbox
{\ifodd\pageno\rightheadline\else\leftheadline\fi} \hbox to
1in{\hfil}\psill} \tenrm
\def\rightheadline{\hfil{\eightrm Nearly Optimal Patchy
Feedbacks }\hfil\folio}
\def\leftheadline{\folio\hfil{\eightrm A. Bressan and F. S. Priuli}\hfil}
\def\onepageout#1{\shipout\vbox
\offinterlineskip \vbox to 2in{\rhead} \vbox to \pageheight
\advancepageno}

\null

\centerline{\bigbf Nearly Optimal Patchy
Feedbacks}
\v
\c{\bigbf for Minimization Problems with Free Terminal Time}
\vskip 2cm \centerline{\it Alberto Bressan$^{(*)}$ and
Fabio S.~Priuli$^{(**)}$}
\v
\c{(*) Department of Mathematics, Penn State University}
\c{University Park, Pa. 16802 U.S.A.}
\c{bressan@math.psu.edu}
\v
\c{(**) Department of Mathematical Sciences, NTNU}
\c{Trondheim, NO-7491, NORWAY}
\c{priuli@math.ntnu.no}
\vsk
\n{\bf Abstract.} The paper is concerned with a general optimization problem
for a nonlinear control system, in the presence of
a running cost and a terminal cost,
with free terminal time.   We prove the existence of
a patchy feedback whose trajectories are all nearly optimal solutions,
with pre-assigned accuracy.
\vsk
\n{\medbf 1 - Introduction}
\v
Consider a general optimization problem
$$
\min_{T,u(\cdot)}~ \left\{ \psi\big(x(T)\big) +
\int_0^T L\big(x(t),\, u(t)\big)\,dt\right\}\,,
\eqno(1.1)
$$
for a nonlinear control system of the form
$$
\dot x = f(x,u)\qquad\qquad u(t)\in \bfU\,.
\eqno(1.2)
$$
Here $x\in\R^n$ describes the state of the system, the upper dot denotes
a derivative w.r.t.~time, and~$\bfU\subset\R^m$ is the set of admissible
control values.  The minimum is sought over all times $T\geq 0$ and all
measurable control functions $u:[0,T]\mapsto \bfU$.

In the literature, several results are
available, which provide  the existence of an optimal control
$t\mapsto u^{ opt}(t)$ in open-loop form [14, 16, 23],
for any fixed initial condition
$$x(0)= y\in\R^n\,.
\eqno(1.3)
$$
On the other hand,  the existence and regularity of an optimal control
in feedback form is a far more difficult issue.
In an ideal situation, one would like to construct a
(sufficiently regular)
feedback $u=U(x)$ such that all trajectories of the corresponding
O.D.E.
$$\dot x= f\big(x, U(x)\big)\eqno(1.4)$$
are optimal w.r.t.~the cost criterion (1.1).
Only few general results are
presently known in this direction [7, 16, 20, 26].
In general, the
optimal feedback can be discontinuous, with
an extremely complicated structure
[8, 17].
Moreover, its  performance may not be robust: an
arbitrarily small external perturbation may produce
trajectories which are far from being optimal [24].

An alternative strategy, pursued in~[3, 18, 19], is to construct
sub-optimal feedbacks, trading off the full optimality in favor of a
simpler structure of the control and the robustness of the resulting
system.
This approach also faces difficulties. In some cases, because of
topological obstructions it is not possible to construct any
continuous asymptotically stabilizing feedback [10, 13, 14, 25], or any
continuous near-optimal feedback [9].
Therefore, one needs to work with discontinuous feedback controls
[11, 12, 21, 22].
For discontinuous O.D.E's, however, no general result
about existence and uniqueness of solutions is available.
Carath\'eodory solutions can be constructed only under
additional assumptions on the structure of discontinuities [15].

Following the approach developed in [1,2,3], asymptotic
stabilization and
optimal control problems can be solved using
{\sl patchy feedbacks} as
discontinuous controls. We recall that a patchy feedback has
a particularly simple structure, since it is a function $u=U(x)$
that is piecewise constant on the state space $\R^n$.
For patchy vector fields, one can prove
that Carath\'eodory solutions forward in time always
exist [1].  Moreover, the set of forward solutions
is stable w.r.t.~small perturbations [2].
The analysis in [3] showed that
any minimum time problem can be approximately solved using these patchy
feedbacks.

Aim of the present paper is to extend the results in [3]
to the general optimization problem (1.1).  In addition, we present a
construction which greatly simplifies the previous approach,
thus clarifying the main lines of the proof.
\v
For convenience, we list here all the
basic assumptions.
\v
\item{\bf (A)} The set of admissible control values
$U\subset\R^m$ is a compact, the function $f:\R^n\times U\mapsto\R^n$ is
continuous w.r.t.~both variables, and twice continuously differentiable
w.r.t.~$x$. In addition, $f$ satisfies the sub-linear
growth condition
$$\big|f(x,u)\big|\leq C_f\big(1+|x|\big)\qquad\qquad\forall\, u\in \bfU\,,
\eqno(1.5)$$
for some constant $C_f\,$.
Both the terminal cost  $\psi:\R^n\mapsto\R$
and the running cost $L:\R^n\times \bfU\mapsto\R\,$
are continuous and non-negative. Moreover, $L$ is strictly positive:
$$
L(x,u)\geq \alpha_0>0\qquad\qquad \forall ~x\in\R^n,~~~u\in\bfU
\,.
\eqno(1.6)$$

\v
Throughout this paper,  $V$ denotes the value function for the
optimization problem (1.1)-(1.2), namely
$$
V(y)\doteq \inf_{T,\,x(\cdot,u)}~ \left\{ \psi\big(x(T)\big) +
\int_0^T L\big(x(t),\, u(t)\big)\,dt\right\}\,,
\eqno(1.7)$$
where the minimization is taken over all $T\geq 0$ and all solutions of
$t\mapsto x(t,u)$, corresponding to a measurable control $u:[0,T]\mapsto\bfU$.
Our main result can be stated as follows.
\v
\n{\bf Theorem 1.}
{\it Let the functions $\psi, L, f$ in (1.1)-(1.2) satisfy the
assumptions {\bf (A)}.
Let $\ve>0$ and a compact set $K\subset\R^n$ be given.
Then there exist a
closed terminal set $S\subseteq\R^n$ and a patchy feedback $u=U(x)$
defined on the
complement $\R^n\setminus S$ such that the following holds.
 For each $y\in K$, every Carath\'eodory solution of
$$\dot x=f\big(x,\,U(x)\big)\,,\qquad\qquad  x(0)=y\eqno(1.8)$$
reaches the set $S$ within finite time.
Calling $\tau\doteq \inf\big\{t\,;~~x(t)\in S\big\}$ the first time where
the trajectory reaches $S$, we have
$$
\psi\big(x(\tau)\big) + \int_0^\tau L\big(x(t),\,
U(x(t))\big)\,dt~\leq  V(y)+\ve\,. \eqno(1.9)
$$
}
\v
We recall that, by well known properties of patchy vector fields, for
every initial point $y\in\R^n\setminus S$ the O.D.E. (1.8) has at
least one forward Carath\'eodory solution.
According to (1.9), all of the solutions starting from the compact set
$K$ are nearly optimal, for the cost (1.1).
\v
In the remainder of the paper, Section 2 contains a brief review of
the main definitions and properties of patchy feedbacks and
patchy vector fields.
The proof of Theorem 1 is then worked out in Section 3.
\vsk
\n{\medbf 2 - Review of patchy feedbacks}
\v
The following definitions were introduced in [1].
\v
\n{\bf Definition 1.} By a {\bf patch} we mean a pair
$\big(\Omega,\, g\big)$ where
$\Omega\subset\R^n$
is an open domain with smooth boundary
$\partial\Om,$ and $g$ is a Lipschitz continuous vector field
defined on a neighborhood of the closure
$\ov\Omega$ of $\Om,$ which points strictly inward at each
boundary point $x\in\partial\Om$.
\v
Calling $\bfn(x)$ the outer normal at the boundary
point $x$, and denoting the inner product by a dot, we thus require
$$
\bfn(x)\cdot g(x)<0\qquad\forall x\in\partial\Omega.
\eqno(2.1)
$$
\v
\n{\bf Definition 2.}  \ We say that $g:\Omega\mapsto\R^n$
is a {\bf patchy vector field} on the open domain $\Om$
if there exists a family of patches
$\big\{ (\Omega_\alpha,~g_\alpha) ;~~ \alpha\in\A\big\}$ such that

\n - $\A$ is a totally ordered set of indices,

\n - the open sets $\Omega_\alpha$ form a locally finite covering of
$\Omega$,

\n - the vector field $g$ can be written in the form
$$
g(x) = g_\alpha (x)\qquad \hbox{if}\qquad x \in
\Omega_\alpha \setminus
\displaystyle{\bigcup_{\beta > \alpha} \Omega_\beta}.
\eqno(2.2)
$$
We shall occasionally adopt the longer notation
$\big(\Omega,\ g,\ (\Omega_\alpha,\,g_\alpha)_{_{\alpha\in \A}} \big)$
to indicate a patchy vector field, specifying both the domain
and the single patches.
\v
By setting
$$
\alpha^*(x) \doteq \max\big\{\alpha \in \A~;~~ x \in \Omega_\alpha
\big\},
\eqno(2.3)
$$
we can write (2.2) in the equivalent form
$$
g(x) = g_{_{\alpha^*(x)}}(x) \qquad \forall  t~~x \in \Omega.
\eqno(2.4)
$$
\v
\n{\bf Remark 1.} It is important to observe
that the patches $(\Omega_\alpha,\,g_\alpha)$
are not uniquely determined by a patchy vector field~$(\Omega,\, g)$.
Indeed,
whenever $\alpha<\beta$, by (2.2) the values of $g_\alpha$ on the set
$\Omega_\alpha\cap\Omega_\beta$ are irrelevant.  Of course, the values of
$g_\alpha$ for $x$ outside the domain $\Omega$ don't matter either.
Therefore, if the open sets
$\Omega_\alpha$
form a locally finite covering of $\Omega$ and
if for each
$\alpha\in \A$ the vector field $g_\alpha$ satisfies
$$\bfn_\alpha(x)\cdot g_\alpha(x)<0\qquad\qquad \forall
x\in\Omega\cap\partial\Omega_\alpha\setminus \bigcup_{\beta>\alpha}
\Omega_\beta\,,\eqno(2.5)$$
then the vector field $g$
in (2.2) is still a patchy vector field.  Indeed,
without changing the function $g$,
one can suitably redefine
the values of each $g_\alpha$ on the set
$\bigcup_{\beta>\alpha} \Omega_\beta$, or outside $\Omega$,
and achieve the strict inequality
$$\bfn_\alpha(x)\cdot g_\alpha(x)<0\qquad\qquad \forall
x\in\partial\Omega_\alpha\,.$$

\v
\n{\bf Remark 2.}  For convenience, we are always
assuming that the single patches
$\Omega_\alpha$ are open, while the vector fields $g_\alpha$
are defined on the closure $\ov\Omega_\alpha\,$.  In certain situations,
it would be natural to choose patches of the form
$$\Omega_1\doteq \big\{x\in\Omega\,;~~\bfn\cdot x<c\big\}\,,\qquad\qquad
\Omega_2\doteq \big\{x\in\Omega\,;~~\bfn\cdot x>c\big\}\,,$$
for some unit vector $\bfn$.
In this way, however, the union $\Omega_1\cup\Omega_2$
does not cover all of $\Omega$, because it does not contain points
where $\bfn\cdot x=c$.    This situation is easily fixed, replacing
$\Omega_1$
by a slightly larger open set which contains also these boundary points.
The resulting vector field
$$g(x)=\cases{g_1(x)\qquad &if\quad $\bfn\cdot x\leq c\,$,\cr
g_2(x)\qquad &if\quad $\bfn\cdot x>c\,$,\cr}$$
can still be written in patchy form.
\vs
If $g$ is a patchy vector field, the differential equation
$$
\dot x = g(x)
\eqno(2.6)
$$
has several useful properties. There are collected in the following theorem,
proved in [1].
\v
\n{\bf Theorem 2.}
{\it  Let $g$ be a patchy vector field.
Then the set
of Carath\'eodory solutions of
(2.6) is closed (in the topology of uniform convergence) but possibly not
connected.  For each  Carath\'eodory solution $t\mapsto x(t) \,$,
the map  $t \mapsto \a^*(x(t))\,$  defined by (2.3) is left-continuous
and non-decreasing.

Given an initial condition
$$
x(t_0)=x_0\,,
\eqno(2.7)
$$
the Cauchy problem (2.6)-(2.7) has at least one forward
solution and at most one backward solution, in the Carath\'eodory
sense.}
\v

\n{\bf Remark 3.}
In some situations it is convenient to adopt a more general definition
of  patchy vector field than the one formulated above. Indeed,
one can consider patches $(\Omega_\alpha,~g_\alpha)$
where the boundary of the domain $\Omega_\alpha$
is only piecewise smooth.   For example, $\Omega_\alpha$ could be a
polytope, or the intersection between a ball and finitely many
half-spaces.
In this more general case, the inward-pointing condition (2.1)
can be reformulated
by asking that, for each boundary point
$x\in\partial \Omega_\alpha$, the vector
$g_\alpha(x)$ lies in the the interior of the
tangent  cone to  $\Omega_\alpha$  at the point $x$.
Namely
$$
g_\alpha(x)\in
{\rm int}\,T_{\!\Omega_\alpha}(x)\,.
\eqno(2.8)
$$
As in [4],
this tangent cone is
defined by
$$
T_{\Omega_\alpha}(x)\doteq
\bigg\{
v\in \R^n~:~ \liminf_{t \downarrow 0}
{d\big(x+tv,\ \Omega_\alpha\big)\over t}=0
\bigg\}.
$$
One can easily check
that all the results concerning patchy
vector fields stated in Theorem 2
remain valid with this more general formulation.

\v
\n{\bf Definition 3.}  \ Let $\big(\Omega,\ g,\
(\Omega_\alpha,\,g_\alpha)_{_{\alpha\in \A}} \big)$ be a patchy
vector field. Assume that there exist control values
$v_\alpha\in \bfU$ such that, for each $\alpha\in\A,$ there holds
$$
g_\alpha(x) = f(x,\, v_\alpha)\qquad\qquad\forall~  x \in
\Omega_\alpha \setminus \bigcup_{\beta > \alpha} \Omega_\beta\,.
\eqno(2.9)
$$
Then the piecewise constant map
$$
U(x) \doteq  v_\alpha\qquad \hbox{if}\qquad x
\in\Omega_\alpha \setminus \bigcup_{\beta > \alpha} \Omega_\beta\,.
\eqno(2.10)
$$
is called a  {\bf patchy feedback}  control on $\Omega\,$.
\v
Recalling (2.3), the patchy feedback control can thus be written on the form
$$U(x) = v_{\alpha^*(x)}\,,$$
\vsk
\n{\medbf 3 - Proof of the theorem}
\v
The proof of Theorem 1 will be given in several steps.
\v
\n{\bf 1.}
Various reductions can be performed.
By a smooth approximation,
we can assume that $\psi\in\C^\infty$.
Moreover, approximating the cost function
$L$ by a more regular function, it is not restrictive to assume that $L$
is twice
continuously differentiable w.r.t.~$x$.  Recalling that
$L(x,u)\geq \alpha_0>0$, we can now replace
$f(x,u)$ by
$$g(x,u)~ \doteq~
{f(x,u)\over L(x,u)}\,,\eqno(3.1)$$
and consider the equivalent problem
$$\inf_{\tau, \,u(\cdot)} \Big\{ \tau +\psi\big(x(\tau)\big)\Big\}\,,
\eqno(3.2)
$$
with dynamics
$$\dot x= g(x,u)\,,\qquad\qquad x(0)=y\,.
$$
Notice that the function $g$ in (3.1) is continuous
w.r.t.~both variables $x,u$, and twice
continuously differentiable w.r.t.~$x$.  Moreover it satisfies the
growth condition
$$\big|g(x,u)\big|\leq {C_f\over\alpha_0}\,
\big(1+|x|\big)\qquad\qquad\forall\, u\in \bfU\,.
$$
In the following, we thus assume without loss of generality that
the running cost is simply
$ L(x,u)\equiv 1$, so that the minimization problem (1.1) reduces to (3.2).
\v
\n{\bf 2.}
Choose a constant $M$ such that
$$M\geq 1\,,\qquad\qquad M\geq \max_{x\in K} \,\psi(x)\,. \eqno(3.3)$$

To fix the ideas, throughout the following
we  assume that $0<\ve<1/8$ and that the compact set
$K$ is contained in the open ball
$B_{\rho}$ centered
at the origin with radius $\rho$.
Because of the sub-linear growth condition (1.5), for $\tau\leq 2M$,
every  trajectory of the
system (1.2) starting at a point $y\in K\subset B_\rho$ will satisfy
the a priori bound
$$\big|x(t)\big|~<~\bar \rho\qquad\qquad\hbox{for all}~~
t\in [0,\tau]\,\subseteq \,[0, 2M]
\,,
\eqno(3.4)$$
where
$$
\bar\rho\doteq e^{C_f \cdot 2M}(\rho+1)\,.
\eqno(3.5)$$
\v
\n{\bf 3.} Let $V=V(y)$ be
the value function for the optimization problem
(3.2), with dynamics (1.2).
We claim that $V$ is semi-concave.  More precisely,
there exists a constant $\kappa$ such that,
for any
$y,y'\in B_{\bar\rho}$, one has
$$V(y') \leq V(y)+\bfw\cdot (y'-y) +\kappa\,{|y'-y|^2\over 2}\,.\eqno(3.6)$$
for some vector $\bfw\in D^+V(y)$
in the upper gradient of $V$ at the point $y$.

Indeed, from the theory of optimal control [6]
it is well known that the optimization problem
(3.2), (1.2) with initial data
$x(0)=y$ has at least one solution,
within the class of chattering controls.  Let
$t\mapsto x(t)=x(t\,;~ y, \tilde u, \tilde\theta)$ be
an optimal chattering trajectory,
with
$$x(0)=y\,,\qquad\qquad
\dot x(t)=\sum_{i=0}^n \theta_i(t) \,f\big(x(t), u_i(t)\big)
\qquad\qquad t\in [0,\tau]\,,\eqno(3.7)$$
for some measurable functions $(\tilde u,\tilde\theta)=
(u_0, \ldots, u_n, \theta_0, \ldots,\theta_n)$ satisfying
$$u_i:[0,\tau]\mapsto\bfU\,,\qquad\qquad \theta_i:[0,\tau]\mapsto[0,1]\,,
\qquad\sum_{i=0}^n \theta_i(t)\equiv 1\,.\eqno(3.8)$$
For any other initial data $y'$, we can consider the same chattering
control $(\tilde u,\tilde \theta)$, always stopping at the same
terminal time $t=\tau$.  This yields the cost
$$V^{\tilde u,\tilde\theta,\tau}(y')=\tau+\psi\big(
x(\tau\,;\, y', \tilde u, \tilde\theta)\big)\,.\eqno(3.9)$$
The regularity assumptions on $f,\psi$ w.r.t.~the variable $x$ imply that,
as $y'$ varies in the ball $B_{\bar \rho}\,$,
the map $y'\mapsto V^{\tilde u,\tilde\theta,\tau}(y')$ is twice continuously
differentiable. Moreover, its $\C^2$ norm
remains bounded:
$$\big\|V^{\tilde u,\tilde\theta,\tau}
\big\|_{\C^2(B_{\bar\rho})}\leq \kappa\,.
\eqno(3.10)$$
Since $\tau\in [0,\,T_{max}]$ while both
$\tilde u$ and $\tilde\theta$ in (3.8) range over compact sets,
this bound is uniform, i.e.~in (3.10) we can take a constant $\kappa>1$
which
does not depend on the particular chattering control,
or on the time $\tau$.
Observing that
$$V(y)=V^{\tilde u,\tilde\theta,\tau}(y)\,,\qquad\qquad
V(y')\leq V^{\tilde u,\tilde\theta,\tau}(y')
\qquad\qquad \forall y'\in B_{\bar\rho}\,,$$
the inequality (3.6) follows from (3.10),
choosing $\bfw=\nabla V^{\tilde u,\tilde\theta,\tau}(y)$.
\v
\n{\bf 4.}
As shown in the previous step, the value function
$$V(y)~=~\min_{\tilde u,\tilde \theta,\tau} \,
V^{\tilde u,\tilde \theta,\tau}(y)$$
is Lipschitz continuous on the ball $B_{\bar\rho}$.
In fact, the constant $\kappa>1$ in (3.10) also provides a
Lipschitz constant for $V$, namely
$$V(x)-V(y)\leq \kappa\,|x-y|\qquad\qquad \forall x,y\in B_{\bar\rho}\,.
\eqno(3.11)$$

By Rademacher's theorem, $V$ is differentiable almost everywhere.
At each point  $x\in B_{\bar\rho}$
where the gradient $\nabla V(x)$ exists,
if $V(x)<\psi(x)$ then one has the well known relation [5, 10, 16]
$$\min_{u\in\bfU}~\big\{ \nabla V(x)\cdot f(x,u)\big\}+1=0\,.\eqno(3.12)$$

Consider the open set
$$\D\doteq\big\{ x\,;~~V(x)<\psi(x)\big\}\,.$$
Given $\delta>0$, we can choose finitely many points $y_1,\ldots y_m\in
B_{\bar\rho}\cap\D$
such that $\nabla V(y_i)$ is well defined for each $i=1,\ldots,m$,
and moreover
$$B_{\bar\rho}\cap\D\subseteq \bigcup_{i=1}^m B(y_i, \delta)\,.\eqno(3.13)$$
Define the approximate value function
$$W(x)\doteq \min~\big\{ \psi(x)\,,~ W_1(x)\,,~\ldots~,~ W_m(x)\big\}\,,
\eqno(3.14)$$
where
$$W_i(x)\doteq V(y_i)+\nabla V(y_i)\cdot (x-y_i)+\kappa\,|x-y_i|^2.
\eqno(3.15)$$

We claim that, by
choosing $\delta>0$ sufficiently small, for all
$x\in B_{\bar\rho}$ the following
relations hold.
$$V(x)\leq W(x)\leq V(x)+\ve\,,\eqno(3.16)$$
$$\left|\min_{u\in\bfU}~\big\{ \nabla W_i(x)\cdot f(x,u)\big\}+1
\right|\leq\ve\qquad\qquad \hbox{whenever}~~W_i(x)=W(x)\,,
\eqno(3.17)$$

Indeed, the first inequality in (3.16) follows from
(3.6).  Next, since $f$ is continuous and $\bfU$ is compact,
we can find $\delta_1\in \,]0,1]$
 such that the following conditions hold.
If $x\in B_{\bar\rho}\,$, $\bfw=\nabla V(y)$ exists and
$$\min_{u\in\bfU}~\big\{ \bfw\cdot f(y,u)\big\}+1=0\,,$$
$$|\bfw'-\bfw|\leq 2\kappa\delta_1\,,\qquad\qquad |x-y|\leq\delta_1\,,$$
then
$$\Big|\min_{u\in\bfU}~\big\{ \bfw'\cdot f(x,u)\big\}+1\Big|\leq\ve\,.
\eqno(3.18)$$

We now choose $\delta>0$ such that
$$2\kappa\delta+\kappa\delta^2\leq  \min\Big\{\ve\,,~
{\kappa\,\delta_1^2\over 2}\Big\}\,. $$
Given any $x\in B_{\bar\rho}$, if $j$ is an index such that
$|x-y_j|\leq\delta$, recalling the Lipschitz condition (3.11) we find
$$W(x)~\leq ~V(y_j)+\big|\nabla V(y_j)\big|\,|x-y_j|+
\kappa |x-y_j|^2~\leq~ V(x)+2\kappa |x-y_j|+\kappa|x-y_j|^2\,,$$
$$W(x)-V(x)~\leq~\min~\Big\{ \ve\,,~{\kappa\delta_1^2\over 2}\Big\}\,.
\eqno(3.19)$$
This already yields (3.16).
Comparing  (3.6) with (3.15) we notice that
$$W_i(x)-V(x)\geq \kappa\, {|x-y_i|^2\over 2}\,.$$
Hence from (3.19) it follows
$$|x-y_i|\leq\delta_1\,,
\qquad\qquad \hbox{whenever}~~W_i(x)=W(x)\,.\eqno(3.20)$$
Observing that, if $W(x)=W_i(x)$,
$$\big|\nabla W_i(x)-\nabla W_i(y_i)\big|~\leq~ 2\kappa \,|x-y_i|~\leq~
2\kappa\,\delta_1\,,$$
from (3.18) we deduce the inequality (3.17). This establishes our claim.
\v
\n{\bf 5.}
By the definition of $W_i$, it is clear that all level sets
where $W_i$ is constant are spheres.  Indeed, for any given
constant $c$ we can write
$$\big\{ x\,;~~W_i(x)=c\big\}~=~
\big\{x\,;~~|x-x_i|=r\big\}\,,$$
with $x_i=y_i-\nabla V(y_i)/2\kappa$ and a suitable radius $r$.

For each $i=1,\ldots,m$, consider the set
$$\D_i\doteq \big\{x\in B_{\bar\rho}\,;~
~W_i(x)=W(x)\big\}\,.\eqno(3.21)$$
In this step we show that there exists a minimum
radius $r_{min}>0$ and a maximum radius $r_{max}$
such that, fixed $x\in\D_i$,
the level set where $W_i=W_i(x)$ is a sphere of center $x_i$ and radius $r$ with
$$0<r_{min}\leq r\leq r_{max}\,.\eqno(3.22)$$
Indeed, since $\ve<1/2$, by (3.17) it follows
$$\big|\nabla W_i(x)\big|\,\big|f(x,u)\big|>{1\over 2}\,.\eqno(3.23)$$
Calling
$$M_f\doteq\max_{|x|\leq\bar\rho, u\in\bfU} \big|f(x,u)\big|\leq
C_f\,(1+\bar\rho)\,,$$
from (3.23) we deduce
$$\big|\nabla W_i(x)\big|> {1\over 2M_f}\,.$$
Therefore, for any $\xi$ such that $W_i(\xi)=W_i(x)$,
$$|\xi-x_i|=|x-x_i|~=~{\big|\nabla W_i(x)\big|\over 2\kappa}~>~
{1\over 4\kappa\,M_f}~\doteq~ r_{min}\,.$$

On the other hand, by (3.15) and (3.20) we have
$$\big|\nabla W_i(x)\big|\leq  \big|\nabla W_i(y_i)\big|+2\kappa
|x-y_i|~\leq~\kappa+2\kappa\delta_1~\leq~ 3\kappa\,.$$
Hence, for any $\xi$ such that $W_i(\xi)=W_i(x)$,
$$|\xi-x_i|=|x-x_i|~=~{\big|\nabla W_i(x)\big|\over 2\kappa}~\leq {3\over 2}~
\doteq~ r_{max}\,.\eqno(3.24)$$
\v
\n{\bf 6.}
We are now ready to construct the near-optimal
patchy feedback. We will define $U(x)$ on the open set
$$\Omega\doteq\big\{x\in B_{\bar\rho}\,;~~W(x)<\psi(x)\big\}\,,
\eqno(3.25)$$
and the required terminal set $S$ will be defined as
$S\doteq\R^n\setminus\Omega$.
Given $\eta>0$ small, for each point $x\in \D_i$ consider the
point (see Figure~1)
$$p_i^x~\doteq ~ {2\over 3}\, x + {1\over 3}\, x_i+\eta{x-x_i\over |x-x_i|}
$$
and the ball $B_i^x=B\big(p_i^x,~|x-x_i|/3\big)$ centered at $p_i^x$
with radius $r=|x-x_i|/3\,$.
By (3.17), there exists a nearly-optimal control value $u=u_i^x\in \bfU$
such that
$$\nabla W_i(x)\cdot f(x,u_i^x)\leq -1+\ve\,.\eqno(3.26)$$
Consider the lens-shaped region
$$\Gamma_i^x\doteq B\Big(p_i^x,~{|x-x_i|\over 3}\Big)
\setminus \ov B\big(x_i\,,~ |x-x_i|-\eta\big).
\eqno(3.27)$$
Its upper boundary will be denoted as
$$\partial^+\Gamma_i^x\doteq \partial\Gamma_i^x\setminus \ov B(x_i\,,~
|x-x_i|-\eta\big)\,.\eqno(3.28)$$
Moreover, for $z\in \partial^+\Gamma_i^x$, we write
$\bfn_i(z)$ for the outer unit normal at the point $z$.

\midinsert
\vskip 10pt
\centerline{\hbox{\psfig{figure=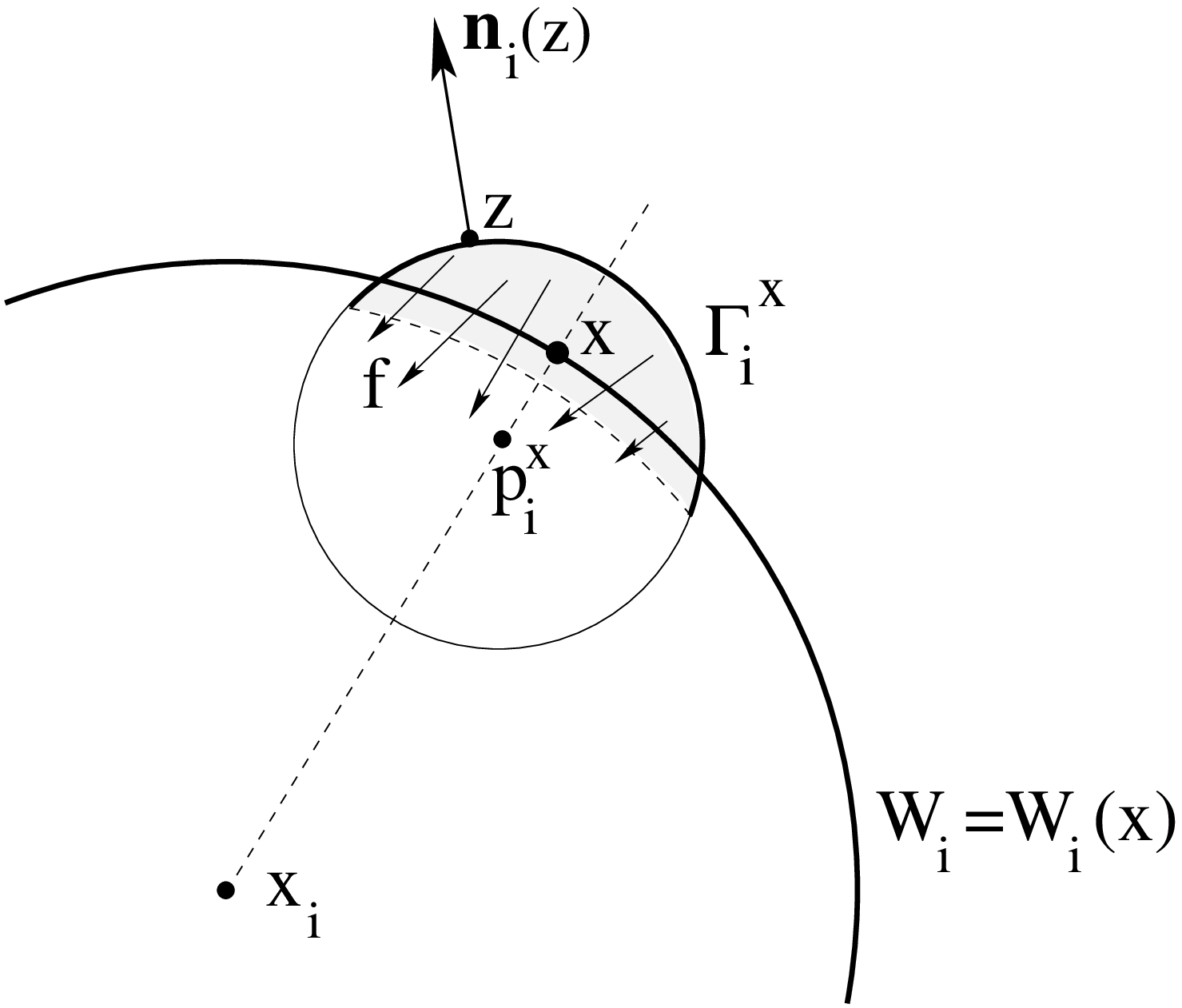,width=6cm}}}
\c{Figure 1. Construction of a lens-shaped patch.}
\vskip 10pt
\endinsert

We claim that, by choosing $\eta>0$ sufficiently small, the following
holds:
$$
\nabla W_i(z)\cdot f(z, u_i^x)\leq -1+2\ve\qquad\qquad
\forall z\in \Gamma_i^x\,,\eqno(3.29)$$
$$\bfn_i(z)\cdot f(z, u_i^x)\leq -\eta\qquad\qquad
\forall z\in \partial^+\Gamma_i^x\,.\eqno(3.30)$$
Moreover, the constant $\eta>0$ can be chosen
uniformly valid for all $i=1,\ldots,m$ and all $x\in \D_i\,$.

For fixed $i,x$ this is clear because, as $\eta\to 0$,
the diameter of the set $\Gamma_i^x$ approaches zero. Moreover,
as $z$ varies on the upper boundary $\partial^+\Gamma_i^x\,$,
all the unit normals $\bfn_i(z)$  approach the vector $\nabla W_i(x)/
\big|\nabla W_i(x)\big|\,$.   Therefore, both inequalities (3.29)-(3.30)
follow from (3.26).

We now observe that $f=f(x,u)$ is uniformly continuous
on the compact domain $B_{\bar\rho}\times\bfU$. Moreover,
on each set $\D_i$,
the gradient $\nabla W_i(x)$ is uniformly Lipschitz continuous
and bounded away from zero.
Finally,
the radius of each level set, where $W_i$ is constant, by (3.22)
is uniformly bounded above and below.
This allows us to choose a constant $\eta>0$ uniformly valid for all $i,x$.
\v
\n{\bf 7.}
To achieve a nearly optimal feedback, we would need the inequality
$$
\nabla W(z)\cdot f(z, u_i^x)\leq -1+4\ve\qquad\qquad
\forall z\in \Gamma_i^x\,.\eqno(3.31)$$
If $W(z)=W_i(z)$ for all $z\in \Gamma_i^x$,
this is a trivial consequence of (3.29).
However, we must also consider the case where
some of the points $z\in \Gamma_i^x$ lie in a region where
$W(z)=W_j(z)<W_i(z)$, for some different index $j$.
For this purpose, we observe that the set where $W_i=W_j$ is always
a hyperplane, say
$$\H_{ij}\doteq
\big\{ x\,;~~W_i(x)=W_j(x)\big\}=\big\{ x\,;~~\bfn_{ij}\cdot x = c_{ij}\big\}\,.
\eqno(3.32)$$
for a suitable constant  $c_{ij}$ and a unit normal vector $\bfn_{ij}\,$.
The orientation of $\bfn_{ij}$ will be chosen so that
$$\big\{ x\,;~~W_i(x)<W_j(x)\big\}=\big\{ x\,;~~\bfn_{ij}\cdot x < c_{ij}\big\}\,.
$$
We claim that, by choosing $\eta>0$ sufficiently small, uniformly w.r.t.~$i,x$,
one of the following two cases occurs (see Figure~2).
\v
\n CASE 1: At every point $z\in \Gamma_i^x\cap\H_{ij}$ one has
$$\bfn_{ij}\cdot f(z, u_i^x)~<~-\eta\,.\eqno(3.33)$$
\v
\n CASE 2: At every point $z\in \Gamma_i^x$ one has
$$
\nabla W_j(z)\cdot f(z, u_i^x)~\leq~ -1+4\ve\,.\eqno(3.34)$$
\v
Indeed, assume that (3.33) fails. Then there exists a point $z^*\in
\Gamma_i^x\cap \H_{ij}$ such that
$$\bfn_{ij}\cdot f(z^*, u_i^x)~\geq ~-\eta\,.\eqno(3.35)$$
By (3.32) and the orientation of the unit
vector $\bfn_{ij}\,$, we can write
$$\nabla W_j(z^*)=\nabla W_i(z^*)-\beta\,\bfn_{ij}\eqno(3.36)$$
for some constant $\beta>0$.
Together, (3.29) and (3.35) now imply
$$\eqalign{\nabla W_j(z^*)\cdot f(z^*, u_i^x)&=
\nabla W_i(z^*)\cdot f(z^*, u_i^x)-\beta\,\bfn_{ij}\cdot f(z^*, u_i^x)
\cr
&\leq -1+2\ve +\beta\,\eta~\leq~-1+3\ve\,,\cr}\eqno(3.37)$$
provided that we choose $\eta>0$ sufficiently small.
Since $f$ and $\nabla W_j$ are uniformly Lipschitz continuous,
from (3.37) it follows that (3.34) is valid for all $z$ sufficiently close to
$z^*$.  By reducing the size of $\eta>0$, we can make
the diameter of the lens-shaped domain $\Gamma_i^x$ as  small as we like.
Hence the inequality (3.34) will hold for all $z\in \Gamma_i^x\,$.

To prove our claim, it remains to observe that the functions $f$ and
$\nabla W_i$ are uniformly continuous, and that the constant
$\beta$ in (3.36) remains uniformly bounded.   Hence the constant $\eta>0$
can be chosen uniformly valid for all $i,j,x\,$.

We now define the smaller domain
$$\Omega_i^x\doteq \Gamma_i^x\setminus \bigcup_{j\in I_i}\big\{
z\in\R^n\,;~~W_j(z)\leq W_i(z)\big\}\eqno(3.38)$$
where $I_i\subset\{1,\ldots,m\}$ is the set of indices $j\not= i$
for which CASE 1 applies.

By the previous analysis, for each $j$ such that $W(z)=W_j(z)$
for some $z\in \Gamma_i^x$, two cases can occur.
If CASE 1 applies, then the
vector field $f(\cdot, u_i^x)$ is strictly inward-pointing
along the portion of the boundary $\partial \Omega_i^x$
where $W_i=W_j$.
On the other hand, if CASE 2 applies, then (3.34) holds on the entire
domain $\Gamma_i^x\,$.

\midinsert
\vskip 10pt
\centerline{\hbox{\psfig{figure=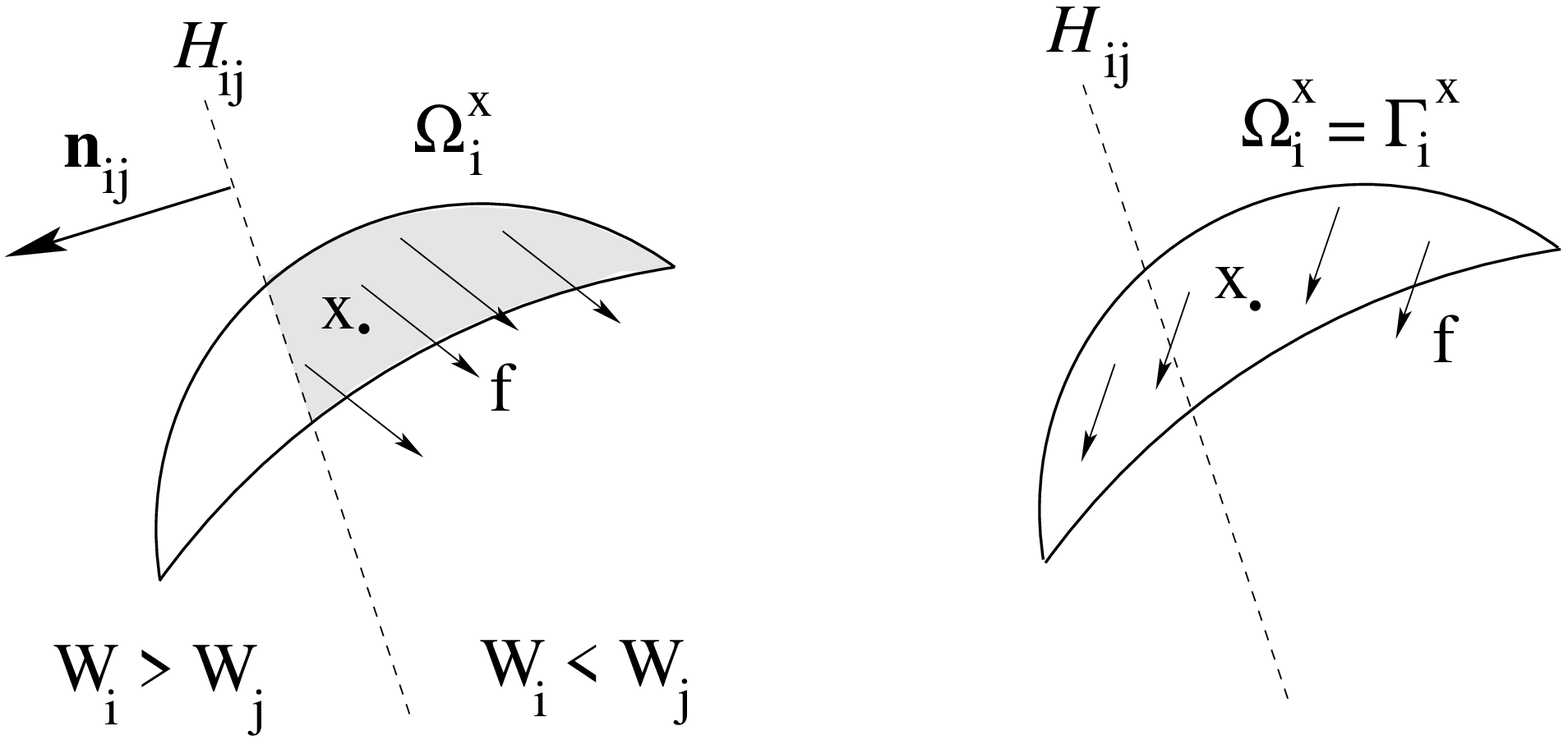,width=10cm}}}
\c{Figure 2. If the domain $\Gamma_i^x$ intersects the half-space
where $W_j<W_i$,
two cases must be considered.}
\c{Left: in Case 1, the vector field $f(\cdot, u_i^x)$
points toward the set where
$W_i<W_j$. As a patch we then}
\c{take the shaded region $\Omega_i^x\subset\Gamma_i^x$.
~Right: In Case 2, $f(\cdot, u_i^x)$
points toward the set where $W_j<W_i$.}
\c{We can now take
$\Omega_i^x=\Gamma_i^x$, because the control $u_i^x$ is
nearly optimal on this whole region.}
\vskip 10pt
\endinsert
\v
\n{\bf 8.}  Consider the family of all domains $\Omega_i^x$, as
$i\in \{1,\ldots,m\}$ and $x$ ranges over the closure of the set
$\Omega\doteq\big\{ x\in B_{\bar\rho}\,;~~W(x)<\psi(x)\big\}$.
%
It now remains to select
finitely many domains $\Omega_i^x$ which cover the compact set $\ov\Omega$.
This last step, however, must be done with some care because on the lower
portion of the boundary
$$\partial^-\Omega_i^x \doteq  \partial
 \Omega_i^x\cap \ov B\big(x_i\,,~ |x-x_i|-\eta\big)\eqno(3.39)$$
the vector field $f(\cdot, u_i^x)$ may not be inward-pointing.
To cope with this problem, we first observe that
there exists a uniform constant $h>0$ such that
$$W_i(z)\leq W_i(x)-h\,,\eqno(3.40)$$
for every $i,x$ and every $z\in\partial^-\Omega_i^x\,$.

We now set $M^*\doteq \max~\big\{ W(x)\,;~x\in B_{\bar\rho}\big\}$,
and split the domain $\Omega$ in sub-domains of the form
$$\Omega_\ell\doteq \big\{ x\in \Omega\,;~~M^*-(\ell+1)h<W(x)<M^* -\ell h
\big\}\,.\eqno(3.41)$$
For each $\ell$, we cover the compact set $\ov\Omega_\ell$
with  finitely many domains $\Omega_i^x$, constructed as
in step {\bf 7}, choosing
$x\in\ov \Omega_\ell\,$.
After a relabelling of both the domains and the correspondent vector
fields from (3.26), this yields the patches (see Figure~3)
$$\big(\Omega_{\ell,\alpha}\,,~f(\cdot, u_{\ell,\alpha})\big)\,,\qquad\qquad
\alpha=1,\ldots, N_\ell\,.\eqno(3.42)$$
On the collection of all patches (3.42) we define the lexicographic order:
$$(\ell,\alpha)\prec (\ell',\alpha')\qquad
\hbox{iff\qquad either $\ell<\ell'$~~
or~~$\ell=\ell'$~and $\alpha<\alpha'$}.$$

\midinsert
\vskip 10pt
\centerline{\hbox{\psfig{figure=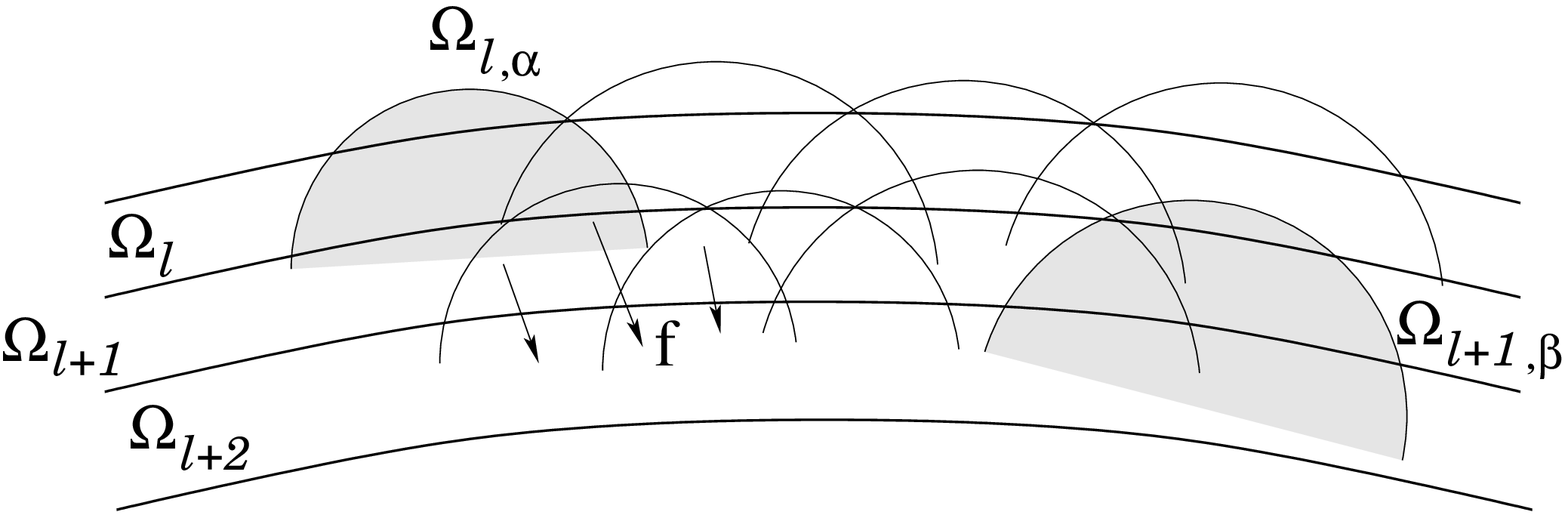,width=10cm}}}
\c{Figure 3. The domain $\Omega=\cup\,\Omega_\ell$ is covered
by a family of patches $\Omega_{\ell,\alpha}\,$,
ordered like tiles on a roof.}
\vskip 10pt
\endinsert

We claim that the above construction yields a patchy vector field:
$$g(x)\doteq f(x,u_{\ell,\alpha})\qquad\hbox{iff}\qquad
x\in \Omega_{\ell,\alpha}\setminus\bigcup_{(\ell,\alpha)\prec (\ell',\alpha')}
\Omega_{\ell',\alpha'}\eqno(3.43)$$
Indeed, according to Remark 1,
it suffices to check that, for each patch $\Omega_{\ell,\alpha}=\Omega_i^x$,
the vector field $f(\cdot, u_{\ell,\alpha})=f(\cdot, u_i^x)$
is inward pointing at every
point of the set
$$\Omega\cap\partial\Omega_{\ell,\alpha}
\setminus\bigcup_{(\ell,\alpha)\prec (\ell',\alpha')}
\Omega_{\ell',\alpha'}\,.$$
In the present case, this is clear, because the only boundary points
where $f(\cdot, u_i^x)$ is not inward pointing are
those on the lower boundary $\partial^-\Omega_i^x$.
Since $x\in \ov \Omega_\ell$ we have $W(x)\leq M^*-\ell h$,
and hence by (3.40)
$$W(z)\leq M^*-(\ell+1)h\qquad\qquad  \forall
z\in \partial^-\Omega_i^x\,.$$
Therefore, given any point $z\in\partial^-\Omega_i^x\,$, either
$W(z)=\psi(z)$ and $z\notin\Omega$, or else
$z$ is
contained in a patch $\Omega_{\ell',\alpha'}$ with $\ell'>\ell$,
as required in Remark 1.
\v
\n{\bf 9.}
To complete the proof, we now check that the patchy
feedback that we have constructed is nearly optimal.
We recall that, by the analysis in step {\bf 7}, for every $i,x$ we have
$$\nabla W_i(z)\cdot f(z, \, u_i^x)\leq -1+4\ve
\qquad\qquad \forall z\in \Omega_i^x\,.\eqno(3.44)$$

Now
take any initial point $y\in K$ and let $t\mapsto x(t)$ be any
Carath\'eodory solution of the Cauchy problem
$$\dot x=g(x)\qquad\qquad x(0)=y\,,$$
with $g$ defined at (3.43).
If $y\in K\setminus\Omega$, $\Omega= \big\{x\in
B_{\bar\rho}\,;~~W(x)<\psi(x)\big\}\,$ as in (3.25), we are in the
terminal set $S$ and there will be no
evolution, since it is more convenient to stay in $y$ than to move
along a trajectory.
Otherwise, call $\tau\geq 0$ the first time at which $x(t)$ reaches the boundary
of the set $\Omega$. By (3.44) we have
$$W\big(x(\tau)\big)-W(y)~=~\int_0^\tau
\Big[{d\over dt} W\big(x(t)\big)\Big]\,dt
~\leq ~(-1+4\ve)\tau\,,$$
hence
$$\tau~\leq~ {W(y)-W\big(x(\tau)\big)\over 1-4\ve}~\leq~ 2 W(y)~\leq~
2\psi(y)~\leq 2M\,.$$
By (3.4),
it follows that $x(\tau)$ cannot be on the boundary of $B_{\bar\rho}$.
We thus conclude that
$W\big(x(\tau)\big)=\psi\big(x(\tau)\big)$.
Stopping at time $\tau$, since $W(x)\geq 0$
and $V(y)\leq M$, the total cost can be
estimated as
$$\tau +\psi\big(x(\tau)\big)~\leq~ {W(y)-W\big(x(\tau)\big)\over 1-4\ve}
+\psi\big(x(\tau)\big)~\leq~
{V(y)+\ve\over 1-4\ve}~\leq V(y)+ {\ve(M+1)\over 1-\ve}\,.$$
Since $\ve>0$ was arbitrary, this completes the proof.

\vsk
\centerline{\medbf References}
\v
\i{[1]}
F.~Ancona and  A.~Bressan,
   Patchy vector fields and asymptotic stabilization.
  {\it ESAIM - Control, Optimiz. Calc. Var.
  } {\bf 4} (1999), 445--471
\v
\i{[2]}
F.~Ancona and  A.~Bressan,
   Flow stability of patchy vector fields and
  robust feedback stabilization
  {\it SIAM J. Control Optim.
  } {\bf 41} (2003), 1455--1476.
\v
\i{[3]}
F.~Ancona and  A.~Bressan,
   Nearly time optimal stabilizing patchy feedbacks.
  {\it Ann. Inst. Henri Poincare, Analyse Non Lineaire},
{\bf 24} (2007), 279-310.
\v
\i{[4]}
J.~P. Aubin and A.~Cellina,
{\it  Differential inclusions. Set-valued maps and viability theory.}
Springer-Verlag, Berlin, 1984.
\v
\i{[5]} M.~Bardi and I.~Capuzzo Dolcetta,
{\it Optimal Control and Viscosity solutions of Hamilton-Jacobi-Bellman
Equations}. Birkh\"auser, Boston, 1997.
\v
\i{[6]} L.~D.~Berkovitz, {\it Optimal Control Theory}.
Springer-Verlag, New York, 1974.
\v
\i{[7]} V.~G.~Boltyanskii,
    Optimal feedback controls
  {\it SIAM J. Control Optim.
  } {\bf 4} (1966), 326--361.
\v
\i{[8]} U.~Boscain and B.~Piccoli,
{\it Optimal Syntheses for Control Systems on 2-D Manifolds}.
Springer-Verlag, Berlin, 2004.
  \v
\i{[9]} A.~Bressan,
    Singularities of stabilizing feedbacks,
  {\it Rend. Sem. Mat. Univ. Politec. Torino,
    } {\bf 56} (1998), 87--104.
\v
\i{[10]} R.~W.~Brockett,
      Asymptotic stability and feedback stabilization,
in {\it Differential Geometric Control Theory},
R.W. Brockett, R.S. Millman, and H.J. Sussmann
Eds., 1983. Birkha\"user Boston, 1983.
\v
\i{[11]}  F.~H.~Clarke, Yu.~S.~Ledyaev, L.~Rifford and R.~J.~Stern,
Feedback stabilization and Lyapunov functions
{\it SIAM J. Control Optim.} {\bf 39} (2000), 25--48.
\v
\i{[12]} F.~H.~Clarke, Yu.~S.~Ledyaev, R.~J.~Stern and P.~R.~Wolenski,
{\it Nonsmooth Analysis and Control Theory} Springer-Verlag
New York, 1998.
\v
\i{[13]}
J.~M.~Coron,
A necessary condition for feedback stabilization.
{\it Systems Control Lett.} {\bf  14}  (1990), 227--232.
\v
\i{[14]}
A.~F.~Filippov,  On certain questions in the theory of optimal control.
{\it J. SIAM Control} {\bf  1},  (1962), 76-84.
\v
\i{[15]} A.~F.~Filippov, {\it Differential Equations with Discontinuous
Righthand Sides}. Kluwer, Dordrecht, 1988.
\v
\i{[16]} W.~Fleming and R.~Rishel, {\it
Deterministic and Stochastic Optimal Control}.
Springer-Verlag, Berlin, 1975.
\v
\i{[17]} I.~Kupka,
The ubiquity of Fuller's phenomenon.
In:  {\it Nonlinear Controllability and Optimal Control},
H.~J.~Sussmann Ed.,
In: {\it Nonlinear Controllability and Optimal Control}, pp.~313--350,
Marcel Dekker, New York, 1990.
\v
\i{[18]} H.~Ishii and S.~Koike,
   On $\epsilon$-optimal controls for state constraints problems
  {\it  Ann. Inst. H. Poincar\'e - Analyse Non Lin\'eaire
  } {\bf 17} (2000),
  473--502.
\v
\i{[19]}  S.~Nobakhtian and R.~J.~Stern,
   Universal near-optimal feedbacks.
  {\it   J. Optim. Theory Appl.
  } {\bf 107} (2000), 89--122
  \v
\i{[20]}
    B.~Piccoli and H.~Sussmann,
   Regular synthesis and sufficiency conditions for optimality.
  {\it  SIAM J. Control Optim.
  } {\bf 39}, (2000), 359--410.
  \v
\i{[21]}  L.~Rifford,
      Stratified semiconcave control-Lyapunov functions and the
stabilization problem.
{\it  Ann. Inst. H. Poincar\'e - Analyse Non lin\'eaire} {\bf 22} (2005),
343-384.
  \v
\i{[22]}  J.D.~Rowland and R.~B.~Vinter,
    Construction of optimal feedback controls.
  {\it   Systems Control Lett.
  } {\bf 16} (1991), 357--367.
  \v
\i{[23]} E.~D.~Sontag,
{\it Mathematical Control Theory. Deterministic Finite Dimensional Systems.
         Second edition.} Springer-Verlag,
New York, 1998.
  \v
\i{[24]} E.~D.~Sontag,
   Stability and stabilization: discontinuities and the
  effect of disturbances.
  In {\it  Nonlinear Analysis, Differential Equations, and Control.}
 F.H. Clarke and R.J. Stern Eds., Kluwer, 1999, pp.~551--598.
    \v
\i{[25]} E.~D.~Sontag and  H.~J.~Sussmann, Remarks on continuous feedback,
in; Proc. IEEE Conf. Decision and Control, Albuquerque, IEEE Publications,
Piscataway 1980, pp.~916-921.
\v
\i{[26]} H.~J.~Sussmann,
Synthesis, presynthesis, sufficient conditions for optimality and
subanalytic sets.
In:  {\it Nonlinear Controllability and Optimal Control},
H.~J.~Sussmann Ed.,
Marcel Dekker,  New York, 1990, pp.~1--19.

\bye